\documentclass[12pt]{article}
\usepackage{indentfirst}

\newtheorem{theorem}{Theorem}
\newtheorem{proposition}[theorem]{Proposition}
\newtheorem{lemma}[theorem]{Lemma}
\newtheorem{corollary}[theorem]{Corollary}

\def\be{\begin{equation}}
\def\ee{\end{equation}}
\def\eqref#1{(\ref{#1})}
\def\diag{{\rm diag}}
\def\Uniform{{\rm Uniform}}

\def\newpartitle#1{\bigbreak 
  \noindent {\bf #1.} \ \ 
  \nopagebreak}
\def\proof{\newpartitle{Proof}}
\def\remark{\newpartitle{Remark}}

\def\ackn{\newpartitle{Acknowledgements}}

\def\darkbox{{\vrule height9pt width7pt depth1pt}}
\def\qed {\hfill\darkbox\\[-7pt]}
\def\eqqed{\eqno\darkbox}

\def\trace{\mbox{trace}}
\def\NR{\vspace{5pt}\cr}

\def\V{{\cal V}}
\def\P{{\bf P}}
\def\E{{\bf E}}

\def\IR{{\bf R}}
\def\Var{{\rm Var}}
\def\Cov{{\rm Cov}}
\def\half{{1 \over 2}}
\def\un{\underbar}
\def\bone{{\bf 1}}
\def\bzero{{\bf 0}}
\def\inn#1#2{ \langle #1, \, #2 \rangle }

\def\overN{{1 \over N}}
\def\fhat{\widehat{f}}

\begin{document}

\baselineskip=18pt

\begin{center}\Large\bf Efficiency of Reversible MCMC Methods:\ \ Elementary 
 \\[4pt]
 Derivations and Applications to Composite Methods
\end{center}
\bigskip
\centerline{\large 
             \ Radford M.\ Neal \ \ \ \ \ \ \ \ \ \ Je{f}frey S.\ Rosenthal}
\centerline{\small\texttt{radford@utstat.utoronto.ca}\hspace{24pt}
                  \texttt{jeff@math.toronto.edu~~~}}
\medskip
\centerline{\sl Department of Statistical Sciences, University of Toronto}
\smallskip
\bigskip
\centerline{(May, 2023; revised March 2024)}
\bigskip
\medskip

\begin{quotation}\noindent
\textbf{Abstract.}\ \ We review criteria for comparing the efficiency
of Markov chain Monte Carlo (MCMC) methods with respect to the
asymptotic variance of estimates of expectations of functions of
state, and show how such criteria can justify ways of combining
improvements to MCMC methods.  We say that a chain on a finite state
space with transition matrix $P$ efficiency-dominates one with
transition matrix $Q$ if for every function of state it has lower (or
equal) asymptotic variance.  We give elementary proofs of some
previous results regarding efficiency dominance, leading to a
self-contained demonstration that a reversible chain with transition
matrix $P$ efficiency-dominates a reversible chain with transition
matrix $Q$ if and only if none of the eigenvalues of $Q-P$ are
negative.  This allows us to conclude that modifying a reversible MCMC
method to improve its efficiency will also improve the efficiency of a
method that randomly chooses either this or some other reversible method,
and to conclude that improving the efficiency of a reversible update
for one component of state (as in Gibbs sampling) will improve the
overall efficiency of a reversible method that combines this and other updates.
It also explains how antithetic MCMC can be more efficient
than i.i.d.\ sampling.
We also establish conditions that can guarantee that a method
is not efficiency-dominated by any other method.
\end{quotation}

\smallskip

\section{Introduction}

Markov chain Monte Carlo (MCMC) algorithms
(e.g.~\cite{mcmchandbook})
estimate the expected
value of a function $f:S\to\IR$ with respect to a probability
distribution $\pi$ on a state space $S$, which in this paper we assume
to be finite, using an estimator such as
$$
\fhat_N \ = \ \overN \sum_{k=1}^N f(X_k)
\, ,
$$
where $X_1,X_2,X_3,\ldots$ is a time-homogeneous Markov chain with stationary 
distribution $\pi$, having transition probabilities $P(x,y)$ from state $x$
to state $y$ (often viewed as a matrix $P$).  

An important measure of the efficiency of this
estimator is its \un{asymptotic variance}:
\be\label{vdef}
\!\!\!v(f,P)
\ := \ \lim_{N\to\infty} N \, \Var\Big[ \fhat_N \Big]
\ = \ \lim_{N\to\infty} N \, \Var\Big[ \overN \sum_{i=1}^N f(X_i) \Big]
\ = \ \lim_{N\to\infty} \overN \, \Var\Big[ \sum_{i=1}^N f(X_i) \Big]
\, .\ \ \ 
\ee
For the irreducible Markov chains used for MCMC, the initial state
of the chain does not affect the asymptotic variance, and the bias
of the estimator converges to zero at rate $1/N$ regardless of initial state. 
(In practice, an initial portion of the chain is usually simulated but
not used for estimation, in order to reduce the bias in a
finite-length run.)

If we run the chain for a large number of iterations, $N$, we therefore
expect that $v(f,P)/N$ will be an indication of the likely squared error 
of the estimate obtained. Indeed, when $v(f,P)$ is finite, one can show
(e.g., \cite[Theorem 5]{tierney}) that a Central Limit Theorem applies, with 
the distribution of $(\fhat-E_{\pi}(f))\,/\,\sqrt{v(f,P)/N}$ converging to 
$N(0,1)$.

We are therefore motivated to try to modify the chain to reduce
$v(f,P)$, ideally for all functions $f$.
We say that one transition matrix, $P$,
\un{efficiency-dominates} another one, $Q$, if
$$
v(f,P) \ \le \ v(f,Q)
\qquad \textrm{for \ all}
\ f:S\to\IR
\, .
$$
Various conditions are known
\cite{tierney, geyerstatsci, radfordnonrev}
which ensure that $P$ efficiency-dominates $Q$.
One of these, for reversible chains, is the \un{Peskun-dominance} condition
\cite{peskun, tierney2}
% that $P(x,S) \ge Q(x,S)$ whenever $x\not\in S$,
which on a finite state space is that $P(x,y) \ge Q(x,y)$ for all $x\not=y$.
This condition is widely cited and has gotten significant
recent attention \cite{LiEtAl,GagnonMaire,Zanella},
and even extended to non-reversible chains \cite{AndrieuLivingstone}.
But it is a very strong condition, and
$P$ might well efficiency-dominate $Q$ even if it does
not Peskun-dominate it.

In this paper, we focus on reversible chains with a finite state
space.  We present several known equivalences of efficiency dominance,
whose proofs were previously scattered in the literature, 
sometimes only hinted at, and sometimes based on very technical
mathematical arguments.  We provide complete \un{elementary} proofs of
them in Sections~\ref{sec-varform}, \ref{sec-equiv},
and~\ref{sec-effequivproof}, using little more than simple linear
algebra techniques.  

In Section~\ref{sec-compareresults}, we use these equivalences to
derive new results, which can show efficiency dominance for
some chains constructed by composing multiple component transition
matrices, as is done for the Gibbs Sampler.
These results are applied to methods
for improving Gibbs sampling in a companion paper~\cite{radfordnew}.
In Section~\ref{sec-eigenresults}, we consider eigenvalue connections,
and show how one can
sometimes prove that a reversible chain cannot be efficiency-dominated
by any other reversible chain, and
also explain (Corollary~\ref{antitheticcor})
how antithetic MCMC can be more efficient than i.i.d.\ sampling, 
These results also allow an easy re-derivation,
in Section~\ref{sec-peskun}, of the fact that Peskun dominance implies
efficiency dominance.

\section{Background Preliminaries}
\label{sec-prelim}

We assume that the state space $S$ is finite, with $|S|=n$, and
let $\pi$ be a probability distribution on $S$,
with $\pi(x)>0$ for all $x\in S$, and $\sum_{x\in S}\! \pi(x) = 1$.
For functionals \mbox{$f,g:S\to\IR$}, define the $L^2(\pi)$ inner product by
$$
\inn{f}{g} \ = \ \sum_{x\in S} f(x) \, g(x) \, \pi(x) \, .
$$
That is, $\inn{f}{g} = \E_\pi[f(X) \, g(X)]$.
Equivalently, if we let $S=\{1,2,\ldots,n\}$,
represent a function $f$ by the column vector
$f=\big[ f(1),\ldots,f(n) \big]^T$,
and let $D=\diag(\pi)$ be the $n \times n$ diagonal matrix with
$\pi(1),\,\ldots,\pi(n)$ on the diagonal,
then $\inn{f}{g}$ equals the matrix product~$f^T\! D g$.

We aim to estimate expectations with respect to $\pi$ by using a
time-homogeneous Markov chain $X_1,X_2,X_3,\ldots$ on $S$,
with transition probabilities 
\mbox{$P(x,y) = \P(X_{t+1}\!=\!y \, | \, X_t\!=\!x)$,}
often written as a matrix $P$, for which $\pi$ is a \un{stationary
distribution} (or \un{invariant distribution}):
$$
 \pi(y)\ =\ \sum_{x\in S} \pi(x) P(x,y)
$$
Usually, $\pi$ is the only stationary distribution, though we
sometimes consider transition matrices that are not irreducible 
(see below), for which this is not true, as building-blocks for other chains.

For $f:S\to\IR$, let $(Pf):S\to\IR$ be the function defined by
$$
(Pf)(x)
\ = \ \sum_{y\in S} P(x,y) \, f(y) \, .
$$
Equivalently, if we represent $f$ as a vector of its values
for elements of $S$, then $Pf$ is the product of the matrix $P$ 
with the vector $f$.  Another interpretation is that 
$(Pf)(x) = E_P[f(X_{t+1})|X_t=x]$, where $E_P$ is expectation 
with respect to the transitions defined by $P$.
We can see that
$$
\inn{f}{Pg} \ = \ \sum_{x\in S} \sum_{y\in S} f(x) \, P(x,y) \, g(y) \, \pi(x)
\, .
$$
Equivalently, $\inn{f}{Pg}$ is the matrix product $f^T\! D P g$.
Also, $\inn{f}{Pg} = \E_{\pi,P}[f(X_t) \, g(X_{t+1})]$,
where $\E_{\pi,P}$ means expectation with respect to the Markov chain
with initial state drawn from the stationary distribution $\pi$
and proceeding according to $P$.

A transition matrix $P$ is called \un{reversible} with respect to
$\pi$ if $\pi(x) \, P(x,y) = \pi(y) \, P(y,x)$ for all $x,y\in S$.
This implies that $\pi$ is a stationary distribution for $P$, since
$\sum_x \pi(x) P(x,y) = \sum_x \pi(y) P(y,x) =
\pi(y) \sum_x P(y,x) = \pi(y)$.  

If $P$ is reversible,
$\inn{f}{Pg} = \inn{Pf}{g}$ for all $f$ and $g$ --- i.e., $P$ is
\un{self-adjoint}
(or, \un{Hermitian}) with respect to $\inn{\cdot}{\cdot}$.
Equivalently, $P$ is reversible with respect to $\pi$
if and only if $DP$ is a symmetric matrix ---
i.e., $DP$ is self-adjoint with respect to the classical dot-product.
This allows us to easily verify some well-known facts about
reversible Markov chains:

\begin{lemma}
If $P$ is reversible with respect to $\pi$
then: (a) the eigenvalues of $P$ are real; (b)
these eigenvalues can be associated with real eigenvectors; (c) if
$\lambda_i$ and $\lambda_j$ are eigenvalues of $P$ with
$\lambda_i \ne \lambda_j$, and $v_i$ and $v_j$ are real eigenvectors
associated with $\lambda_i$ and $\lambda_j$, then $v_i^T D v_j = 0$,
where $D$ is the diagonal matrix with $\pi$ on the diagonal
(i.e., $\inn{v_i}{v_j}=0$); (d) all the eigenvalues of $P$ are in $[-1,1]$.
\end{lemma}

\proof
Since $DP$ is symmetric and $D$ is diagonal, $DP=(DP)^T=P^TD$.
(a) If $Pv\,=\,\lambda v$ with $v$
non-zero, then $\overline\lambda\, \overline v^T \,=\, \overline v^T\!
P^T$, hence $\overline \lambda\,(\overline v^T\! D v)\,=\, \overline
v^T\!  P^T D v \,=\, \overline v^T\! DP v\,=\,\lambda\, (\overline
v^T\! D v)$.  Since $\overline v^T\! D v$ is non-zero (because $D$ has
positive diagonal elements), it follows that $\overline
\lambda=\lambda$, and hence $\lambda$ is real. (b) If $Pv\,=\,\lambda
v$, with $P$ and $\lambda$ real and $v$ non-zero, then at least one of
$\mbox{Re}(v)$ and $\mbox{Im}(v)$ is non-zero and is a real
eigenvector associated with~$\lambda$. (c) $\lambda_i(v_i^T D
v_j)\,=\, v_i^T P^T D v_j\, =\, v_i^T DP v_j\, =\, \lambda_j (v_i^T D
v_j)$, which when $\lambda_i\ne\lambda_j$ implies that $v_i^T D
v_j=0$. (d) Since rows of $P$ are non-negative and sum to one, the
absolute value of an element of the vector $Pv$ can be no larger than
the largest absolute value of an element of $v$.  If $Pv=\lambda v$,
this implies that $|\lambda| \le 1$, hence $\lambda\in[-1,1]$.
\qed

The self-adjoint property
implies that $P$ is a ``normal operator'', which guarantees
(e.g., \cite[Theorem~2.5.3]{horn}) the existence of an
\un{orthonormal basis}, $v_1,v_2,\ldots,v_n$, of eigenvectors for $P$,
with $Pv_i = \lambda_i v_i$ for each $i$, and $\inn{v_i}{v_j}=\delta_{ij}$.
(In particular, this property implies that $P$ is \un{diagonalisable}
or \un{non-defective}, but it is stronger than that.)
Without loss of generality, we can take $\lambda_1=1$,
and $v_1 = \bone := \big[ 1,1,\ldots,1 \big]^T$,
so that $v_1(x) = \bone(x) = 1$ for all $x\in S$,
since $P \bone = \bone$ due to the transition probabilities in $P$
summing to one.
We can assume for convenience that all of $P$'s eigenvalues
(counting multiplicity) satisfy
$\lambda_1 \ge \lambda_2 \ge \ldots \ge \lambda_n$.

In terms of orthonormal eigenvectors of $P$, any functions $f,g:S\to\IR$
can be written as linear combinations $f = \sum_{i=1}^n a_i v_i$
and $g = \sum_{j=1}^n b_j v_j$.
It then follows from orthonormality of these eigenvectors that
$$ \inn{f}{g} = \sum_i a_i b_i,\ \ \
   \inn{f}{f} = \sum_i (a_i)^2,\ \ \ 
   \inn{f}{Pg} = \sum_i a_i b_i \lambda_i,\ \ \
   \inn{f}{Pf} = \sum_i (a_i)^2 \lambda_i.
$$
% NEW:
In particular, $\inn{f}{v_i} = a_i$, so the $a_i$ are the
projections of $f$ on each of the $v_i$. This shows
that $\sum_{i=1}^n v_i v_i^T D$ is equal to the identity matrix,
since for all $f$,
$$ \Big(\sum_{i=1}^n v_i v_i^T D\Big) f \ =\
   \sum_{i=1}^n v_i (v_i^T D f)
\ = \ \sum_{i=1}^n v_i \, \inn{v_i}{f}
\ = \ \sum_{i=1}^n a_i v_i \ =\  f
\, .
$$
Furthermore, any self-adjoint $A$ whose eigenvalues are all zero must be 
the zero operator, since we can write any $f$
% using an orthonormal basis of eigenvectors of $A$
as $f=\sum_{i=1}^n a_i v_i$, from which it follows that
$Af = \sum_{i=1}^n \lambda_i a_i v_i = 0$.

A matrix $A$ that is self-adjoint with respect to $\inn{\cdot}{\cdot}$
has a \un{spectral representation} in
terms of its eigenvalues and eigenvectors as 
$A = \sum_{i=1}^n \lambda_i v_i v_i^T D$. If $h:\IR\to\IR$ we can define
$h(A) := \sum_{i=1}^n h(\lambda_i) v_i v_i^T D$, which is easily seen to
be self-adjoint.  Using $h(\lambda)=1$ gives the identity matrix.
One can also easily show that $h_1(A)+h_2(A)=(h_1+h_2)(A)$
and $h_1(A)h_2(A)=(h_1h_2)(A)$, and hence (when all $\lambda_i\ne0$) that
$A^{-1}=\sum_{i=1}^n \lambda_i^{-1} v_i v_i^T D$ and (when all $\lambda_i\ge0$)
that $A^{1/2}=\sum_{i=1}^n \lambda_i^{1/2} v_i v_i^T D$, so both of these
are self-adjoint.  
Finally, note that if $A$ and $B$ are self-adjoint, so is $ABA$.

We say that $P$ is \un{irreducible} if movement from any $x$ to any $y$
in $S$ is possible via some number of transitions that have positive
probability under $P$.  An irreducible chain will have only one stationary
distribution.  A reversible irreducible $P$ will have
$\lambda_i<1$ for $i\ge2$. (As an aside, this implies that $P$ is
\un{variance bounding}, which in turn
implies that $v(f,P)$ from~\eqref{vdef} must be finite for each $f$
\cite[Theorem~14]{varbound}.)
For MCMC estimation, we want our chain to be irreducible, but
irreducible chains are sometimes built using transition matrices that
are not irreducible --- for example, by letting $P=(1/2)P_1+(1/2)P_2$,
where $P_1$ and/or $P_2$ are not irreducible, but $P$ is irreducible.

An irreducible $P$ is \un{periodic} with period $p$ if $S$ can be partitioned 
into $p>1$ subsets $S_0,\ldots,S_{p-1}$ such that $P(x,y)=0$ if $x\in S_a$ and
$y \notin S_b$, where $b=a\!+\!1\ \mbox{mod}\ p$ (and this is not true for
any smaller $p$). Otherwise, $P$ is
\un{aperiodic}.  An irreducible periodic chain that is reversible must have 
period 2, and will have $\lambda_n=-1$ and $\lambda_i>-1$ for $i \ne n$.
A reversible aperiodic chain will have all $\lambda_i>-1$.

Since $v(f,P)$ as defined in~\eqref{vdef} only involves variance,
we can subtract off the mean of $f$ without affecting the asymptotic variance.
Hence, we can always assume without loss of generality that $\pi(f)=0$,
where $\pi(f) := \E_\pi(f) = \sum_{x\in S} f(x) \, \pi(x)
= \inn{f}{\bone}$.  In other words, we can assume that $f\in L^2_0(\pi) 
:= \{ f: \pi(f)=0, \ \pi(f^2)<\infty\}$,
where the condition that $\pi(f^2) = \E_\pi(f^2)$ be finite is 
automatically satisfied when $S$ is finite, and hence can be ignored.
Also, if $\pi(f) = 0$,
then $\inn{f}{\bone} = \inn{f}{v_1} = 0$, so $f$ is orthogonal to
$v_1$, and hence its coefficient $a_1$ is zero.

Next, note that
$$
\inn{f}{P^kg}
\ = \ \sum_{x\in S} f(x) \, (P^kg)(x) \, \pi(x)
\ = \ \sum_{x\in S} \sum_{y\in S} f(x) P^k(x,y) g(y) \pi(x)
\ = \ \E_{\pi,P}[ f(X_t) \, g(X_{t+k}) ]\,.
$$
where $P^k(x,y)$ is the $k$-step transition probability from $x$ to $y$.
If $f\in L^2_0(\pi)$ (i.e., the mean of $f$ is zero), this is
the \un{covariance} of $f(X_t)$ and $g(X_{t+k})$, when the chain 
is started in stationarity (and hence is the same for all $t$).  We
define the \un{lag-$k$ autocovariance}, $\gamma_k$, as:
$$
\gamma_k
\ := \ \Cov_{\pi,P}[f(X_t),f(X_{t+k})]
\ := \ \E_{\pi,P}[f(X_t) f(X_{t+k})]
\ = \ \inn{f}{P^kf}, \ \ \ \mbox{for $f\in L^2_0(\pi)$}.
$$
If $f = \sum_{i=1}^n a_i v_i$ as above (with $a_1=0$ since the mean of 
$f$ is zero), then using orthonormality of the eigenvectors $v_i$,
$$
\gamma_k
\ = \
\inn{f}{P^kf}
\ = \ \sum_{i=2}^n \sum_{j=2}^n \inn{a_iv_i}{P^k(a_jv_j)}
\ = \ \sum_{i=2}^n \sum_{j=2}^n \inn{a_iv_i}{(\lambda_j)^k a_jv_j)}
\ = \ \sum_{i=2}^n (a_i)^2 (\lambda_i)^k
\, .
$$
In particular, $\gamma_0 = \inn{f}{f} := \|f\|_{L^2(\pi)} = \sum_i (a_i)^2$.
(If the state space $S$ were not finite, we would need to 
require $f \in L^2(\pi)$, but finite variance is guaranteed with a finite
state space.)

One particular example of a transition matrix $P$, useful for
comparative purposes, is $\Pi$, the operator
corresponding to i.i.d.\ sampling from $\pi$.  It is defined by
$\Pi(x,y) = \pi(y)$ for all $x\in S$.  This operator satisfies
$\Pi \bone = \bone$, and $\Pi f = 0$ whenever $\pi(f)=0$.  Hence,
its eigenvalues are $\lambda_1=1$ and $\lambda_i=0$ for $i\ne1$.

\section{Relating Asymptotic Variance to Eigenvalues}
\label{sec-varform}

In this section, we consider some expressions for the asymptotic
variance, $v(f,P)$, of~\eqref{vdef}, beginning with a result relating
the asymptotic variance to the eigenvalues of $P$.  This result (as
observed by~\cite{geyerstatsci}) can be obtained (at least in the
aperiodic case) as a special case of the more technical results of Kipnis
and Varadhan \cite[eqn~(1.1)]{kipnis}.

\begin{proposition}\label{vevalform}
If $P$ is an irreducible (but possibly periodic) Markov chain
on a finite state space $S$, which is reversible with respect to $\pi$,
with orthonormal basis $v_1,v_2,\ldots,v_n$ of eigenvectors,
and corresponding eigenvalues $\lambda_1\ge\lambda_2\ge\cdots\ge\lambda_n$,
and $f\in L^2_0(\pi)$ with $f=\sum_i a_i v_i$, then
the limit $v(f,P)$ in~\eqref{vdef} exists, and
$$
v(f,P)
\ = \
\sum_{i=2}^n (a_i)^2
\, + \, 2 \, \sum_{i=2}^n\, (a_i)^2 \, {\lambda_i \over 1-\lambda_i}
\ = \
\sum_{i=2}^n\, (a_i)^2 \, {1+\lambda_i \over 1-\lambda_i}
\, .
$$
\end{proposition}

\proof
First, by expanding the square, using stationarity, and
collecting like terms, we obtain the well-known result 
that for $f\in L^2_0(\pi)$,
\begin{eqnarray*}
\overN \, \Var\,\bigg(\sum_{i=1}^N f(X_i) \bigg)
& = &
\overN \, \E_{\pi,P}\left[ \bigg( \sum_{i=1}^N f(X_i) \bigg)^{\!2}\, \right]
\\[3pt]
& = &
\overN \, \bigg( N \, \E_{\pi,P}[f(X_j)^2]
\,+\, 2 \sum_{k=1}^{N-1} (N\!-\!k)\ \E_{\pi,P}[f(X_j) \ f(X_{j+k})]\, \bigg)
  \\[3pt]
% \ = \
% \E_{\pi,P}[f(X_j)^2]
% + 2 \sum_{k=1}^{N-1} {N-k \over N} \ \E_{\pi,P}[f(X_j) \ f(X_{j+k})]
& = & \gamma_0 \,+\, 2 \sum_{k=1}^{N-1} {N-k \over N} \ \gamma_k
\, ,
\end{eqnarray*}
where $\gamma_k = \Cov_{\pi,P}[f(X_j),f(X_{j+k})] = \inn{f}{P^kf}$
is the lag-$k$ autocovariance in stationarity.

Now, $f = \sum_{i=1}^n a_i v_i$, with $a_1=0$ since $\pi(f)=0$, so
$\gamma_k \ = \ \inn{f}{P^kf} \ = \ \sum_{i=2}^n (a_i)^2 (\lambda_i)^k$
and $\gamma_0 \ = \ \sum_{i=2}^n (a_i)^2$. The above then gives that
\be\label{firstvareqn}
\overN \, \Var\left( \sum_{i=1}^N f(X_i) \right)
\ = \
\sum_{i=2}^n (a_i)^2 + 2 \sum_{k=1}^{N-1} {N-k \over N}
\ \sum_{i=2}^n (a_i)^2 (\lambda_i)^k
\, ,
\ee
i.e.\
$$
\overN \, \Var\left( \sum_{i=1}^N f(X_i) \right)
\ = \
\sum_{i=2}^n (a_i)^2 + 2 \sum_{k=1}^\infty \sum_{i=2}^n
I_{k \le N-1} \ {N-k \over N} \ (a_i)^2 (\lambda_i)^k
\, .
$$
If $P$ is aperiodic, then
$\Lambda := \max_{i \ge 2} |\lambda_i| \, < \, 1$,
hence
$\sum_{k=1}^\infty \Big| \sum_{i=2}^n
I_{k \le N-1} \, {N-k \over N} \ (a_i)^2 (\lambda_i)^k \Big|
\le
\sum_{k=1}^\infty \sum_{i=2}^n
\Big| I_{k \le N-1} \, {N-k \over N} \ (a_i)^2 (\lambda_i)^k \Big|
\,\le\, \sum_{k=1}^\infty \sum_{i=2}^n (a_i)^2 (\Lambda)^k
\,=\, \gamma_0 \, \Lambda/(1-\Lambda) < \infty$,
so the above sum is \un{absolutely summable}.
This lets us exchange the limit and summations to obtain
$$
v(f,P) \ := \ \lim_{N\to\infty}
\overN \, \Var\left( \sum_{i=1}^N f(X_i) \right)
\ = \ 
\sum_{i=2}^n (a_i)^2
+ 2 \sum_{i=2}^n \sum_{k=1}^\infty
\lim_{N\to\infty}
\left[ I_{k \le N-1} \ {N-k \over N} \ (a_i)^2 (\lambda_i)^k \right]
$$
$$
\ \ \ \ \ \ = \ 
\sum_{i=2}^n (a_i)^2 + 2 \sum_{i=2}^n \sum_{k=1}^\infty
(a_i)^2 (\lambda_i)^k
\ = \ 
\sum_{i=2}^n (a_i)^2 + 2 \sum_{i=2}^n\,
(a_i)^2 {\lambda_i \over 1-\lambda_i}
\ = \ 
\sum_{i=2}^n\,
(a_i)^2 \, {1+\lambda_i \over 1-\lambda_i}
\, .
$$

If $P$ is periodic, with $\lambda_n=-1$, then
the above $\Lambda=1$, and $\sum_{k=1}^\infty (\lambda_n)^k$
is not even defined, so the above argument does not apply.
Instead, separate out the $i=n$ term in~\eqref{firstvareqn} to get
$$
\overN \, \Var\left( \sum_{i=1}^N f(X_i) \right)
\ = \
\sum_{i=2}^n (a_i)^2 \, + \, 2 \sum_{k=1}^{N-1} {N-k \over N}
\ \sum_{i=2}^{n-1} (a_i)^2 (\lambda_i)^k
\, + \, 2
\sum_{k=1}^{N-1} {N-k \over N} (a_n)^2  (-1)^k
\, .
$$
Since $\Gamma := \max\{|\lambda_2|,|\lambda_3|,\ldots,|\lambda_{n-1}|\} < 1$,
the previous argument applies to the middle double-sum term to show that
$$
\lim_{N\to\infty} 2 \sum_{k=1}^{N-1} {N-k \over N}
\ \sum_{i=2}^{n-1} (a_i)^2 (\lambda_i)^k
\ = \ 
2 \sum_{i=2}^{n-1}\, (a_i)^2\, {\lambda_i \over 1-\lambda_i}
\, .
$$
As for the final term, writing values for $k$ as $2m\!-\!1$ or $2m$, 
we have
\def\IevenN{I_{N \rm\ is\ even}}
\def\IoddNm1{I_{N-1 \rm\ is\ odd}}
$$
\sum_{k=1}^{N-1} {N-k \over N} (-1)^k
\,\ =\, \ 
\overN \sum_{m=1}^{\lfloor (N-1)/2 \rfloor} \Big[-(N-2m+1)+(N-2m)\,\Big]
\ - \ \overN\, \IoddNm1
$$
$$
\,\ = \,\ \overN \sum_{m=1}^{\lfloor (N-1)/2 \rfloor} \Big[-1\Big]
\ - \ \overN\, \IoddNm1
\,\ =\,\ - {\lfloor (N-1)/2 \rfloor \over N}
\ - \ \overN\, \IoddNm1
\, ,
$$
which converges as $N\to\infty$ to
$-\half = {-1 \over 1-(-1)} = {\lambda_n \over 1-\lambda_n}$.
So, we again obtain that
$$
v(f,P)\ =\ \lim_{N\to\infty}
\overN \, \Var\left( \sum_{i=1}^N f(X_i) \right)
\ =\, \
\sum_{i=2}^n (a_i)^2
\, + \, 2 \, \sum_{i=2}^{n-1}\, (a_i)^2 {\lambda_i \over 1-\lambda_i}
\, + \, 2\, (a_n)^2 {\lambda_n \over 1-\lambda_n}\ \
$$
$$
\ =\, \
\sum_{i=2}^n (a_i)^2
\, + \, 2 \, \sum_{i=2}^n (a_i)^2 {\lambda_i \over 1-\lambda_i}
\ =\, \
\sum_{i=2}^n\, (a_i)^2\, {1+\lambda_i \over 1-\lambda_i}
\, .
\eqqed
$$
\smallskip

Note that when $P$ is periodic, $\lambda_n$ will be $-1$, and the
final term in the expression for $v(f,P)$ will be zero.  Such a
periodic $P$ will have
zero asymptotic variance when estimating the expectation of a function
$f$ for which $a_n$ is the only non-zero coefficient.

When $P$ is aperiodic, we can obtain from
Proposition~\ref{vevalform} the more familiar 
\cite{bartlett,billingsley,changeyer,olle,huang,tierney} 
expression for $v(f,P)$ in terms of sums of
autocovariances, though it is not needed for this paper (and
actually still holds without the reversibility condition 
\cite[Theorem 20.1]{billingsley}):

\begin{proposition}\label{onelagproposition}
If $P$ is a reversible, irreducible, aperiodic Markov chain
on a finite state space $S$
with stationary distribution $\pi$,
and $f\in L^2_0(\pi)$, then
\be\label{gammasum}
v(f,P)\ =\ \lim_{N\to\infty}
{1 \over N} \, \Var_\pi\left( \sum_{i=1}^{N} f(X_i) \right)
\ = \
\gamma_0 + 2 \sum_{k=1}^\infty \gamma_k
\, ,
\ee
where $\gamma_k = \Cov_{\pi,P}[f(X_t),f(X_{t+k})]$
is the lag-$k$ autocovariance of $f$ in stationarity.
\end{proposition}

\proof
Since $\gamma_k = \inn{f}{P^kf}$,
and as above in the aperiodic case
$\Lambda := \sup_{i \ge 2} |\lambda_i| < 1$,
the double-sum is again absolutely summable,
and we compute directly that if $f=\sum_i a_i v_i$ then
$$
\gamma_0 + 2 \sum_{k=1}^\infty \gamma_k
\ = \ \inn{f}{f} + 2 \sum_{k=1}^\infty \inn{f}{P^kf}
\ = \ \sum_i (a_i)^2 + 2 \sum_{k=1}^\infty \sum_i (a_i)^2 (\lambda_i)^k
$$
$$
\ = \ \sum_i (a_i)^2 + 2 \sum_i (a_i)^2 \, \sum_{k=1}^\infty (\lambda_i)^k
\ = \ \sum_i (a_i)^2 + 2 \sum_i (a_i)^2 \, {\lambda_i \over 1-\lambda_i}
\, ,
$$
so the result follows from Proposition~\ref{vevalform}.
\qed

Proposition~\ref{vevalform} gives the following
formula for $v(f,P)$ (see also~\cite[Lemma 3.2]{mirageyer}):

\begin{proposition}\label{vprop}
The asymptotic variance, $v(f,P)$, for the functional $f \in L^2_0(\pi)$
using an \un{irreducible} Markov chain $P$ which is \un{reversible}
with respect to $\pi$ satisfies the equation
$$
v(f,P)
\ = \ \inn{f}{f}\ +\ 2\, \inn{f}{P(I\!-\!P)^{-1} f}
\, ,
$$
which we can also write as 
$v(f,P) = \inn{f}{f}\ +\ 2\, \inn{f}{{P \over I-P} \, f}$,
or as $v(f,P) = \inn{f}{{I+P \over I-P} \, f}$.
\end{proposition}

\proof
Let $P$ have an orthonormal basis $v_1,v_2,\ldots,v_n$, with
eigenvalues $\lambda_1,\lambda_2,\ldots,\lambda_n$,
and using this basis let $f=\sum_i a_i v_i$.  Note that $a_1=0$,
since $f$ has mean zero, so we can ignore $v_1$ and
$\lambda_1$. Define 
$h(\lambda) := \lambda (1-\lambda)^{-1}$.  As discussed in
Section~\ref{sec-prelim}, applying $h$ to the eigenvalues
of $P$ will produce another self-adjoint matrix, with the
same eigenvectors, which will equal $P(I-P)^{-1}$.  Using this, we can write
\begin{eqnarray*}
\inn{f}{f}\ +\ 2\, \inn{f}{P(I\!-\!P)^{-1} f}
& = & \sum_i (a_i)^2
    \ +\ 2 \sum_i \sum_j\, \inn{a_i v_i}{P(I\!-\!P)^{-1} (a_j v_j)}
\\[3pt]
& = & \sum_i (a_i)^2 
    \ +\ 2 \sum_i \sum_j\, \inn{a_i v_i}{\lambda_j (1-\lambda_j)^{-1} (a_j v_j)}
\\[3pt]
& = & \sum_i (a_i)^2\ +\ 2 \sum_i\, (a_i)^2\, \lambda_i\, (1-\lambda_i)^{-1}\, ,
\end{eqnarray*}
so the result follows from Proposition~\ref{vevalform}.
\qed

\remark If we write $P = \sum_{i=1}^n \lambda_i v_i v_i^T D$,
so $I\!-\!P = \sum_{i=1}^n (1-\lambda_i) v_i v_i^T D$, then
on $L^2_0(\pi)$ this becomes $I\!-\!P = \sum_{i=2}^n (1-\lambda_i) v_i v_i^T D$,
so $(I\!-\!P)^{-1} = \sum_{i=2}^n (1-\lambda_i)^{-1} v_i v_i^T D$.

\remark The inverse $(I\!-\!P)^{-1}$ in Proposition~\ref{vprop} is
on the restricted space $L^2_0(\pi)$ of functions $f$ with $\pi(f)=0$.
That is,
$(I\!-\!P)^{-1}\,(I\!-\!P)\,f\,=\,
(I\!-\!P)\,(I\!-\!P)^{-1}\,f\,=\,f$ for any $f$ in $L^2_0(\pi)$.
By contrast, $I\!-\!P$ will \un{not} be invertible
on the full space $L^2(\pi)$ of \un{all} functions on $S$,
since, for example, $(I\!-\!P)\bone=\bone-\bone=\bzero$,
so if $(I\!-\!P)^{-1}$ existed on all of $L^2(\pi)$ then we would have
the contradiction that
$\bone = (I\!-\!P)^{-1} (I\!-\!P) \, \bone
= (I\!-\!P)^{-1} \, \bzero = \bzero$.

\section{Efficiency Dominance Equivalences}
\label{sec-equiv}

Combining Proposition~\ref{vprop} with the definition
of efficiency dominance proves:

\begin{proposition}\label{pinvprop}
For reversible irreducible Markov chain transition matrices $P$ and $Q$,
$P$ efficiency-dominates $Q$ \un{if and only if}
$\inn{f}{P(I\!-\!P)^{-1} f}
\le \inn{f}{Q(I\!-\!Q)^{-1} f}$ for all $f \in L^2_0(\pi)$,
or informally that
$\inn{f}{{P \over I-P}f} \le \inn{f}{{Q \over I-Q}f}$
for all $f \in L^2_0(\pi)$.
\end{proposition}

Next, we need the following fact:

\begin{lemma}\label{effequivlemma}
% If $P$ and $Q$ are self-adjoint $n\times n$ matrices
% with eigenvalues contained in $(-\infty,1)$, then
If $P$ and $Q$ are reversible and irreducible Markov chain transition matrices,
$\inn{f}{P(I\!-\!P)^{-1}f} \le \inn{f}{Q(I\!-\!Q)^{-1}f}$
for all $f \in L^2_0(\pi)$
\un{if and only if}
$\inn{f}{Pf} \le \inn{f}{Qf}$ for all $f \in L^2_0(\pi)$.
\end{lemma}

Lemma~\ref{effequivlemma} follows
from the very technical results of
Bendat and Sherman \cite{bendat}.
It is somewhat subtle since the equivalence is only
for \un{all} $f$ at once, not for individual $f$;
see the discussion after Lemma~\ref{inverselemmaoneway} below.
In Section~\ref{sec-effequivproof} below,
we present an \un{elementary} proof.
(For alternative direct proofs of Lemma~\ref{effequivlemma}
and related facts, see also \cite[Chapter~V]{Bhatia}.)

Combining Lemma~\ref{effequivlemma}
and Proposition~\ref{pinvprop} immediately shows
the following, which is also shown by Mira and 
Geyer~\cite[Theorem 4.2]{mirageyer}:

\begin{theorem}\label{linprop}
For Markov chain transition matrices $P$ and $Q$ that are
reversible and irreducible,
$P$ efficiency-dominates $Q$ \un{if and only if}
$\inn{f}{Pf}
\le \inn{f}{Qf}$ for all $f \in L^2_0(\pi)$,
i.e.\ \un{if and only if}
$\inn{f}{(Q\!-\!P)f} \ge 0$ for all $f \in L^2_0(\pi)$.
\end{theorem}

\remark Here the restriction that $f \in L^2_0(\pi)$, i.e.\
that $\pi(f)=0$, can be omitted, since
if $c := \pi(f) \not= 0$ then $f=f_0+c$ where $\pi(f_0)=0$, and
$\inn{f}{Pf} = \inn{f_0+c}{P(f_0+c)} = \inn{f_0}{Pf_0} + c^2$,
and similarly for $Q$.
But we do not need this fact here.

\remark Some authors (e.g.~\cite{mira}) say that $P$
\un{covariance-dominates} $Q$ if $\inn{f}{Pf} \le \inn{f}{Qf}$ for all
$f\in L^2_0(\pi)$, or equivalently if $\Cov_{\pi,P}[f(X_t), \, f(X_{t+1})]$ is
always smaller under $P$ than under $Q$.  The surprising conclusion of
Theorem~\ref{linprop} is that for reversible chains this is
\un{equivalent} to efficiency dominance --- i.e., to $v(f,P) \le v(f,Q)$
for all $f\in L^2_0(\pi)$.  So, there is no need to consider the two
concepts separately.

\bigskip

To make the condition $\inn{f}{(Q\!-\!P)f} \ge 0$ for all $f$ more concrete,
we have the following:

\begin{lemma}\label{poslemma}
Any self-adjoint matrix $J$ satisfies
$\inn{f}{Jf} \ge 0$ for all $f$
\un{if and only if} the eigenvalues of $J$ are all non-negative,
which is \un{if and only if} the eigenvalues of $DJ$ are all non-negative
where $D=\diag(\pi)$.
\end{lemma}

\proof
Let $J$ have orthonormal basis of eigenvectors $v_1,v_2,\ldots,v_n$ as
in Section~\ref{sec-prelim}, so any $f$ can be written
as $f = \sum_{i=1}^n a_i v_i$.  Then
$$
\inn{f}{Jf}
\ = \ \inn{\sum_i a_i v_i}{\sum_j a_j Jv_j}
\ = \ \inn{\sum_i a_i v_i}{\sum_j a_j \lambda_j v_j}
\ = \ \sum_i (a_i)^2 \lambda_i
\, .
$$
If each $\lambda_i \ge 0$, then this expression must be $\ge 0$.
Conversely, if some $\lambda_i<0$, then choosing $f=v_i$
gives $\inn{f}{Jf} = \lambda_i < 0$.
This proves the first statement.

For the second statement, recall
that $DJ$ is self-adjoint with respect to the classical dot-product.
Hence, by the above, the matrix product $f^TDJf \ge 0$ for all $f$
\un{if and only if} the eigenvalues of $DJ$ are all non-negative.
So, since $f^TDJf = \inn{f}{Jf}$, the two statements are equivalent.
\qed

Combining Lemma~\ref{poslemma} (with $J$ replaced by $Q\!-\!P$) with
Theorem~\ref{linprop} shows:

\begin{theorem}\label{eigenthm}
If $P$ and $Q$ are reversible irreducible Markov chain transitions,
$P$ efficiency-dominates $Q$ \un{if and only if}
the operator $Q-P$ (equivalently, the matrix $Q-P$) has all eigenvalues 
non-negative, which is \un{if and only if} the matrix $D\,(Q-P)$ has all 
eigenvalues non-negative.
\end{theorem}

%\smallskip
\remark By Theorem~\ref{eigenthm}, if $Q-P$ has even a single negative
eigenvalue, say $(Q-P)\,z\, =\, -cz$ where $c>0$,
then there must be some $f\in L^2_0(\pi)$ such that $v(f,Q) < v(f,P)$.
By following through our proof of Lemma~\ref{effequivlemma}
in Section~\ref{sec-effequivproof} below, it might be
possible to construct such an $f$
\un{explicitly} in terms of $z$ and $c$.
We leave this as an open problem.

\smallskip
\remark It might be possible to give another alternative proof of
Theorem~\ref{eigenthm} using the step-wise approach of \cite{radfordnonrev},
by writing $Q\!-\!P = R_1+R_2+\ldots+R_\ell$ where each $R_i$ is of rank one
(e.g., $R_i=\lambda_i v_i v_i^T D$ with $\lambda_i$ and $v_i$ an eigenvalue
and eigenvector of $Q\!-\!P$).
We leave this as another open problem.

\bigskip

Theorem~\ref{eigenthm} allows us to prove the following, 
which helps justify the phrase ``efficiency-dominates''
(see also \cite[Section 4]{mirageyer}):

\begin{theorem}\label{partialord}
Efficiency dominance is a \un{partial order} on reversible chains, i.e.:
\begin{itemize}
\item[(a)]
It is \un{reflexive}: $P$ always efficiency-dominates $P$;
\item[(b)]
It is \un{antisymmetric}:
if $P$ efficiency-dominates $Q$,
and $Q$ efficiency-dominates $P$,
then $P=Q$;
\item[(c)] It is \un{transitive}:
if $P$ efficiency-dominates $Q$,
and $Q$ efficiency-dominates $R$,
then $P$ efficiency-dominates $R$.
\end{itemize}
\end{theorem}

\proof
Statement~(a) is trivial.
Statement~(c) is true because
$v(f,P) \le v(f,Q)$ and $v(f,Q) \le v(f,R)$ imply
$v(f,P) \le v(f,R)$.
For statement~(b), 
Theorem~\ref{eigenthm} implies that both $Q\!-\!P$ and $P\!-\!Q$ have
all eigenvalues non-negative, hence their eigenvalues must all be zero,
which implies (since $Q\!-\!P$ is self-adjoint) that $Q\!-\!P=0$, 
and hence $P=Q$.
\qed

\remark Statement (b) of Theorem~\ref{partialord} does not hold if we
do not assume reversibility.
For example, if $S=\{1,2,3\}$, and $\pi=\Uniform(S)$,
and $P(1,2)=P(2,3)=P(3,1)=1$,
and $Q(1,3)=Q(3,2)=Q(2,1)=1$,
then $v(f,P)=v(f,Q)=0$ for all $f:S\to\IR$,
so they each (weakly) efficiency-dominate the other,
but $P\not=Q$.

% \section{New Efficiency Dominance Results}
% \label{sec-newresults}

\section{Efficiency Dominance of Combined Chains}
\label{sec-compareresults}

Using Theorems~\ref{linprop} and~\ref{eigenthm}, we can now prove some
new results about efficiency dominance that are useful when Markov
chains are constructed by combining two or more chains.

We first consider the situation where we randomly choose to apply
transitions defined either by $P$ or by $Q$.  For example, $P$ might
move about one region of the state space well, while $Q$ moves about a
different region well.  Randomly choosing either $P$ or $Q$ may
produce a chain that moves well over the entire state space.  The
following theorem says that if in this situation we can 
improve $P$ to $P'$, then the random combination will also be
improved:

% \medskip\bf Q1-result. \rm
\begin{theorem}
Let $P$, $P'$, and $Q$ be reversible with respect to $\pi$,
with $P$ and $P'$ irreducible,
and let $0<a<1$.
Then $P'$ efficiency-dominates $P$ \un{if and only if}
$aP'+(1\!-\!a)Q$ efficiency-dominates $aP+(1\!-\!a)Q$.
\end{theorem}

\proof
Since $P$ and $P'$ are irreducible, so are
$aP'+(1\!-\!a)Q$ and $aP+(1\!-\!a)Q$.
So by Theorem~\ref{eigenthm},
$P'$ efficiency-dominates $P$
\un{if and only if}
$P-P'$ has all non-negative eigenvalues, which is
clearly \un{if and only if}
$a(P\!-\!P')\, =\, [aP+(1\!-\!a)Q]\, -\, [aP'+(1\!-\!a)Q]$
has all non-negative eigenvalues, which is
\un{if and only if}
$aP'+(1\!-\!a)Q$ efficiency-dominates $aP+(1\!-\!a)Q$.
\qed

The next result applies to Markov chains built using component
transition matrices that are not necessarily irreducible, such as
single-variable updates in a random-scan Gibbs sampler, again showing
that improving one of the components will improve the combination,
assuming the combination is irreducible:

\begin{theorem}\label{gibbsthm}
Let $P_1,\ldots,P_{\ell}$ and $P'_1,\ldots,P'_{\ell}$ be reversible with
respect to $\pi$ (though not necessarily irreducible).
Let $a_1,\ldots,a_{\ell}$ be mixing probabilities, with $a_k>0$ 
and $\sum_k a_k = 1$, and let
$P = a_1 P_1 + ... + a_{\ell} P_{\ell}$
and $P' = a_1 P'_1 + ... + a_{\ell} P'_{\ell}$.
Then if $P$ and $P'$ are irreducible,
and for each $k$ the eigenvalues of $P_k-P'_k$
(or of $D(P_k-P'_k)$ where $D=\diag(\pi)$)
are all non-negative, then $P'$ efficiency-dominates $P$.
\end{theorem}

\proof
Choose any $f\in L^2_0(\pi)$.
Since $P_k$ and $P'_k$ are self-adjoint,
we have from Lemma~\ref{poslemma} that
% the assumptions imply that
$\inn{f}{(P_k-P'_k)f} \ge 0$
for each $k$.
% or equivalently the matrix product $f^T D (P_k-P'_k) f \ge 0$.
Then, by linearity,
$$
\inn{f}{(P-P')f}
\ = \ \inn{f}{\sum_k a_k \, (P_k-P'_k)f}
\ = \ \sum_k a_k \, \inn{f}{(P_k-P'_k)f}
\ \ge \ 0
\, ,
$$
too.  Hence,
by Theorem~\ref{linprop}, $P'$ efficiency-dominates $P$.
\qed

In the Gibbs sampling application, the state is composed of $\ell$
components, so that $S = S_1 \times S_2 \times \cdots \times S_{\ell}$,
and $P_k$ is the transition that samples a value for component $k$,
independent of its current value, from its conditional distribution
given the values of other components, while leaving the values of these other
components unchanged.  Since it leaves other components unchanged,
such a $P_k$ will not be irreducible.  $P_k$ will be a block-diagonal matrix,
in a suitable ordering of states (different for each $k$), with
$B=|S|/|S_k|$ blocks, each of size $K=|S_k|$. 

For example, suppose $\ell=2$, $S_1=\{1,2\}$, $S_2=\{1,2,3\}$, and
$\pi(x)=1/9$ except that $\pi((1,2))=4/9$.  With lexicographic ordering,
the Gibbs sampling transition matrix for the second component, $P_2$, will be
$$
  P_2\ =\ 
  \pmatrix{ 1/6 & 4/6 & 1/6 & 0 & 0 & 0 \NR
            1/6 & 4/6 & 1/6 & 0 & 0 & 0 \NR
            1/6 & 4/6 & 1/6 & 0 & 0 & 0 \NR
            0 & 0 & 0       & 1/3 & 1/3 & 1/3 \NR
            0 & 0 & 0       & 1/3 & 1/3 & 1/3 \NR
            0 & 0 & 0       & 1/3 & 1/3 & 1/3}\\[-1pt]
$$

Each block of $P_k$ can be regarded as the $K\times K$ transition
matrix for a Markov chain having $S_k$ as its state space, which is
reversible with respect to the conditional distribution on $S_k$ given
the values for other components associated with this block.  For each
block, the eigenvalues and eigenvectors of this transition matrix give rise to 
corresponding eigenvalues and eigenvectors of $P_k$, after prepending
and appending zeros to the eigenvector according to how many blocks
precede and follow this block.
If the transition matrix for each block is
irreducible, there will be $B$ eigenvalues of $P_k$ equal to one, with
eigenvectors of the form $[ 0, \ldots, 0, 1, \ldots, 1, 0 \ldots,
0]^T$, which are zero except for a series of $K$ ones corresponding to
one of the blocks.

The overall transition matrix when using Gibbs sampling to update a
component randomly chosen with equal probabilities will be $P =
(1/\ell) (P_1 + \cdots + P_{\ell})$.  We can try to improve the
efficiency of $P$ by modifying one or more of the $P_k$.  An
improvement to $P_k$ can take the form of an improvement to one of its
blocks, each of which corresponds to particular values of components
of the state other than component $k$.  With each $P_k$ changed to
$P'_k$, the modified overall transition matrix is $P' = (1/\ell) (P'_1
+ \cdots + P'_{\ell})$.

For the example above, we could try to
improve $P$ by improving $P_2$, with the improvement to $P_2$ taking the
form of an improvement to how the second component is changed when the
first component has the value 1, as follows:
$$
  P'_2\ =\ 
  \pmatrix{   0 &  1  &  0  & 0 & 0 & 0 \NR
            1/4 & 2/4 & 1/4 & 0 & 0 & 0 \NR
              0 &  1  &  0  & 0 & 0 & 0 \NR
            0 & 0 & 0       & 1/3 & 1/3 & 1/3 \NR
            0 & 0 & 0       & 1/3 & 1/3 & 1/3 \NR
            0 & 0 & 0       & 1/3 & 1/3 & 1/3}\\[-1pt]
$$
The change to the $3\times3$ upper-left block still leaves it reversible
with respect to the conditional distribution for the second component 
given the value 1 for the first component (which has probabilities 
of $1/6$, $4/6$, $1/6$), but
introduces an antithetic aspect to the sampling.

If we leave $P_1$ unchanged, so $P'_1=P_1$, 
Theorem~\ref{gibbsthm} can be used to show that the $P'$ built with
this modified $P'_2$ efficiency-dominates $P$
built with the original $P_1$ and $P_2$.  The difference $P_2-P'_2$
will also be block diagonal, and its eigenvalues will be those of
the differences in the individual blocks (which are zero for blocks
that have not been changed).  In the example above, the one block 
in the upper-left that changed has difference:
$$
  \pmatrix{1/6 & 4/6 & 1/6 \NR 1/6 & 4/6 & 1/6 \NR 1/6 & 4/6 & 1/6}
  \ -\ 
  \pmatrix{ 0  &  1  &  0  \NR 1/4 & 2/4 & 1/4 \NR  0  &  1  &  0 }
  \ =\ 
  \pmatrix
     {+2/12 & -4/12 & +2/12\NR -1/12 & +2/12 & -1/12\NR +2/12 & -4/12 & +2/12 }
$$
The eigenvalues of this difference matrix are $1/2$, $0$, and $0$.
The eigenvalues of $P_2-P'_2$ will be these plus three more zeros.  If
$P'_1=P_1$, Theorem~\ref{gibbsthm} then guarantees that $P'$
efficiency-dominates $P$, the original Gibbs sampling chain, since
these eigenvalues are all non-negative.  Note that here one cannot 
show efficiency-dominance using Peskun-dominance, since the change
reduces some off-diagonal transition probabilities.

More generally, suppose a Gibbs sampling chain is changed by modifying one
or more of the blocks of one or more of the $P_k$, with the new blocks
efficiency-dominating the old Gibbs sampling blocks (seen as transition matrices
reversible with respect to the conditional distribution for that
block).  Then by Theorem~\ref{eigenthm}, the eigenvalues of the differences
between the old and new blocks are all non-negative, which implies that
the eigenvalues of $P_k-P'_k$ all non-negative for each $k$, which
by Theorem~\ref{gibbsthm} implies that the modified chain efficiency-dominates
the original Gibbs sampling chain.  The practical applications of this
are developed further in the companion paper~\cite{radfordnew}.

\section{Efficiency Dominance and Eigenvalues}
\label{sec-eigenresults}

We will now present some results relating eigenvalues of transition
matrices to efficiency dominance, which can sometimes be used to
show that a reversible transition matrix \un{cannot} be efficiency-dominated
by any other reversible transition matrix.

Say that $P$ \un{eigen-dominates} $Q$ if both $P$ and $Q$ are reversible
and the eigenvalues of $P$ are no greater than the corresponding eigenvalues
of $Q$ --- that is, when the eigenvalues of $P$ (counting multiplicities)
are written non-increasing as
$\lambda_1 \ge \lambda_2 \ge \ldots \ge \lambda_n$,
and the eigenvalues of $Q$ are written non-increasing as
$\beta_1 \ge \beta_2 \ge \ldots \ge \beta_n$,
then $\lambda_i \le \beta_i$ for each $i$.
Then we have (see also \cite[Theorem 3.3]{mirageyer}):

\begin{proposition}\label{eigenefprop}
If $P$ and $Q$ are irreducible and reversible with respect to $\pi$,
and $P$ efficiency-dominates  $Q$, then $P$ eigen-dominates $Q$.
\end{proposition}

\proof
By Theorem~\ref{linprop},
$\inn{f}{Pf} \le \inn{f}{Qf}$ for all $f\in L^2_0(\pi)$.
Hence, the result follows from the
``min-max'' characterisation of eigenvalues
(e.g.\ \cite[Theorem~4.2.6]{horn})
that
$$
\lambda_i \ = \ \inf_{g_1,\ldots,g_{i-1}}
\sup_{f\in L^2_0(\pi) \atop
	{\inn{f}{f}=1 \atop \inn{f}{g_j}=0 \, \forall \, j}}
% \sup_{\substack{f:S\to\IR \\ \|f\|=1 \\ \inn{f}{g_j}=0 \, \forall \, j}
% {\inn{f}{Pf} \over \inn{f}{f}}
\inn{f}{Pf}
\, .
$$
Intuitively, $g_1,\ldots,g_{i-1}$ represent the first $i\!-\!1$ eigenvectors
(excluding the eigenvector~$\bone$ associated with the eigenvalue~1,
since $\bone \notin L^2_0(\pi)$),
so that the new eigenvector $f$ will be orthogonal to them.  However,
since the formula is stated in terms of \un{any} vectors $g_1,\ldots,g_{i-1}$,
the same formula applies for both $P$ and $Q$, thus giving the result.
\qed

The \un{converse} of Proposition~\ref{eigenefprop} does \un{not} hold,
contrary to a claim in \cite[Theorem~2]{mira}.
For example, suppose the state space is
$S=\{1,2,3\}$.  Let $e=0.05$, and let
% define $P$, $Q$, and $R$ as
$$
P\ =\ \pmatrix{\,\half\, & \half   & 0 \NR 
               \,\half\, & \half-e & e \NR 
               \,0\,     & e       & 1\!-\!e\,},\ \ \ \
Q\ =\ \pmatrix{1\!-\!e   & e       & 0 \NR 
               e         & \half-e & \half\, \NR 
               0         & \half   & \half\,},\ \ \ \
R\ =\ \pmatrix{1\!-\!e   & e       & 0 \NR 
               e         & \half   & \half-e\, \NR 
               0         & \half-e & \half+e\,}\vspace{3pt}
$$
\iffalse
$Q(1,2)=Q(2,1)=e$,
$Q(1,1)=1-e$,
$Q(2,2)=\half-e$,
$Q(2,3)=Q(3,2)=Q(3,3)=\half$,
and meanwhile
$P(3,2)=P(2,3)=e$,
$P(3,3)=1-e$,
$P(2,2)=\half-e$,
$P(2,1)=P(1,2)=P(1,1)=\half$.
\fi
These all are reversible with respect to
$\pi=\,$Uniform$(S)$, and are irreducible and aperiodic.
One can see that $P$ eigen-dominates $Q$ (and vice versa),
since $P$ and $Q$ are equivalent upon swapping states~1 and~3, and so
have the same eigenvalues, which are equal (to four decimal places)
to $1$, $0.9270$, $-0.0270$.
However, $P$ does \un{not} efficiency-dominate $Q$,
since $Q\!-\!P$ has eigenvalues $0.7794$, $0$, $-0.7794$ which are
\un{not} all non-negative. (Nor does $Q$ efficiency dominate $P$,
analogously.)

Intuitively, in this example,
$Q$ moves easily between states 2 and 3, but only infrequently
to or from state 1, while
$P$ moves easily between states 1 and 2 but not to or from state 3.
Hence, if, for example, $f(1)=2$ and $f(2)=1$ and $f(3)=3$ so
that $f(1)=\half[f(2)+f(3)]$,
then $v(f,Q) < v(f,P)$, since
$Q$ moving slowly between $\{1\}$ and $\{2,3\}$ doesn't matter, but
$P$ moving slowly between $\{1,2\}$ and $\{3\}$ \un{does} matter.

$R$ is a slight modification to $Q$ that has two smaller off-diagonal
elements, and hence is Peskun-dominated (and efficiency-dominated) by $Q$.
It's eigenvalues are 1, 0.9272, 0.0728, the later two of which are 
strictly larger than those of $P$, so $P$ eigen-dominates $R$.  But the
eigenvalues of $R-P$ are $0.7865$, 0, $-0.6865$, which are not all 
non-negative,
so $P$ does not efficiency-dominate $R$, despite strictly eigen-dominating it.

However, the next result is in a sense a converse of 
Proposition~\ref{eigenefprop} for the special case where 
\un{all} of the non-trivial eigenvalues for $P$ are smaller
than \un{all} of those for $Q$:

\begin{theorem}\label{minmaxeigen}
Let $P$ and $Q$ be irreducible and
reversible with respect to $\pi$, with eigenvalues
$1 = \lambda_1 \ge \lambda_2 \ge \lambda_3 \ge \ldots \ge \lambda_n$ for $P$
and
$1 = \beta_1 \ge \beta_2 \ge \beta_3 \ge \ldots \ge \beta_n$ for $Q$
(counting multiplicities).
Suppose $\max_{i \ge 2} \lambda_i \le \min_{i \ge 2} \beta_i$,
i.e.\ $\lambda_2 \le \beta_n$,
i.e.\ $\lambda_i \le \beta_j$ for \un{any} $i,j \ge 2$.
Then $P$ efficiency-dominates $Q$.
\end{theorem}

\proof
Let $a$ be the maximum eigenvalue of $P$
restricted to $L^2_0(\pi)$, and
let $b$ be the minimum eigenvalue of $Q$
restricted to $L^2_0(\pi)$, which for both $P$ and $Q$ will
exclude the eigenvalue of 1 associated with $\bone$.
The assumptions imply that $a \le b$.
But the ``min-max'' characterisation of eigenvalues
(described in the proof of
Proposition~\ref{eigenefprop} above), applied to the largest
eigenvalue on $L^2_0(\pi)$ (i.e., excluding the eigenvalue~1), implies that
$$
a \ = \ \sup_{f\in L^2_0(\pi) \atop \inn{f}{f}=1} \inn{f}{Pf}
\, .
$$
Also, since $-b$ is the largest eigenvalue of $-Q$,
$$
b \ = \ - \sup_{f\in L^2_0(\pi) \atop \inn{f}{f}=1} \inn{f}{-Qf}
\ = \ \inf_{f\in L^2_0(\pi) \atop \inn{f}{f}=1} \inn{f}{Qf}
\, .
$$
Since $a \le b$, this implies that $\inn{f}{Pf} \le \inn{f}{Qf}$
for any $f\in L^2_0(\pi)$ with $\inn{f}{f}=1$,
and hence (by linearity) for \un{any} $f \in L^2_0(\pi)$.
It then follows from Theorem~\ref{linprop} that $P$ efficiency-dominates $Q$.
\qed

A chain is called ``antithetic'' (cf.~\cite{greenhan}) if all its
eigenvalues (except $\lambda_1=1$) are non-positive,
with at least one negative.  Our next result shows that
such antithetic samplers always efficiency-dominate i.i.d.\ sampling:

\begin{corollary}\label{antitheticcor}
If $P$ is irreducible and
reversible with respect to $\pi$, and has eigenvalues $\lambda_1=1$
and $\lambda_2,\lambda_3,\ldots,\lambda_n \le 0$,
then $P$ efficiency-dominates $\Pi$ (the operator corresponding
to i.i.d.\ sampling from $\pi$).
\end{corollary}

\proof
By assumption, $\max_{i \ge 2} \lambda_i \le 0$.
Also, if
$1 = \beta_1 \ge \beta_2 \ge \beta_3 \ge \ldots \ge \beta_n$ are
the eigenvalues for $\Pi$,
then $\beta_i=0$ for all $i \ge 2$,
so $\min_{i \ge 2} \beta_i = 0$.
Hence, $\max_{i \ge 2} \lambda_i \le \min_{i \ge 2} \beta_i$.
The result then follows from Theorem~\ref{minmaxeigen}.
\qed\vspace{-6pt}

% \proof
% Let $v_1,v_2,\ldots,v_n$ be an orthonormal basis of eigenvectors
% for $P$, with $v_1=\bone$.  Then $P\,\bone=\bone$ and $\Pi\,\bone=\bone$, 
% so $v_1$ is also an eigenvector of $\Pi$ with eigenvalue $\lambda_1=1$.
% Since the rows of $\Pi$ are all equal, the other eigenvalues of $\Pi$ are
% all zero, and can be associated with the same orthonormal eigenvectors 
% $v_2,\ldots,v_n$ as $P$.  So $(\Pi-P)v_1 = 0$, and for
% $i \ge 2$, $Pv_i = \lambda_i v_i$ and $\Pi \, v_i = 0$ so that
% $(\Pi-P)v_i = -\lambda_i v_i$.
% Hence, the eigenvalues of $\Pi-P$ are
% $0,-\lambda_2,-\lambda_2,\ldots,-\lambda_n$, which
% are all non-negative since each $\lambda_i \le 0$.
% Hence, by Theorem~\ref{eigenthm},
% $P$ efficiency-dominates $\Pi$.
% \qed\vspace{-6pt}

\remark
Theorem~1 of \cite{frigessi} shows that if 
$\pi_{\min} = \min_x \pi(x)$, the maximum eigenvalue (other 
than $\lambda_1$) of a transition matrix reversible with respect to
$\pi$ must be greater than or equal to 
$-\pi_{\min}\,/\,(1-\pi_{\min})$, which must be greater than or equal
to $-1/(n\!-\!1)$ since $\pi_{\min} \le 1/n$.
Now, if $f = v_i$ where $\lambda_i \ge -1/(n\!-\!1)$,
then Proposition~\ref{vevalform} gives
$v(f,P) = (1+\lambda_i)/(1-\lambda_i) \ge (n-2)/n$,
which for large $n$ is only slightly smaller than
$v(f,\Pi) = (1+0)/(1-0) = 1$.
% This limits the possible gain from antithetic
% sampling for the worst-case choice of $f$.
On the other hand, we can still have e.g.\ $\lambda_1=1$,
$\lambda_n=-1$, and all the other $\lambda_i=0$,
and then if $f = v_n$ then
$v(f,P) = (1+(-1))/(1-(-1)) = 0/2 = 0$,
which is significantly less than $v(f,\Pi) = 1$.
Hence, the improvement in Corollary~\ref{antitheticcor}
could be large for some functions~$f$, but small for some others.

\bigskip

Since practical interest focuses on whether or not some chain, $P$,
efficiency-dominates another chain, $Q$, Proposition~\ref{eigenefprop}
is perhaps most useful in its contrapositive form --- if $P$
and $Q$ are reversible, and $P$ does \un{not} eigen-dominate $Q$,
then $P$ does \un{not} efficiency-dominate~$Q$.  That is,
if $Q$ has at least one eigenvalue less than the corresponding
eigenvalue of $P$, then $P$ does not efficiency-dominate $Q$.  If
both chains have an eigenvalue less than the corresponding
eigenvalue of the other chain, then neither efficiency-dominates
the other.  

But what if two different chains have exactly the same ordered set
of eigenvalues --- that is, they both eigen-dominate the other?  
In that case, neither efficiency-dominates the other.
To show that, we first prove a result about strict trace comparisons:

% In fact, this remains true even if they just have the same eigenvalue sums:

\begin{theorem}\label{tracesame}
If $P$ and $Q$ are both irreducible transitions matrices, reversible with
respect to $\pi$,
and $P$ efficiency-dominates $Q$,
and $P \ne Q$,
then $\trace(P)<\trace(Q)$, i.e.\ the trace (or equivalently the
sum of eigenvalues) of $P$ is \un{strictly} smaller than that of $Q$.
\end{theorem}

\proof By Theorem~\ref{eigenthm}, if $P$ efficiency-dominates $Q$,
then $Q\!-\!P$ has no negative eigenvalues.  And it cannot have all zero
eigenvalues, since then $Q\!-\!P=0$ (since $Q\!-\!P$ is self-adjoint),
contradicting the premise that $P \ne Q$.  So, $Q\!-\!P$ has at
least one positive eigenvalue, and no negative eigenvalues, and hence
the sum of eigenvalues of $Q\!-\!P$ is strictly positive.  But the sum of
the eigenvalues of a matrix is equal to its trace~\cite[p.~51]{horn}, so 
this is equivalent to $\trace(Q\!-\!P)>0$.  Since trace is linear, this implies
that $\trace(Q)\!-\!\trace(P) > 0$, and hence $\trace(P)<\trace(Q)$.
\qed

\begin{corollary}\label{eigensame}
If $P$ and $Q$ are both irreducible transitions matrices, reversible with
respect to $\pi$, and the eigenvalues
(counting multiplicity)
for both are identical,
and $P \ne Q$, then 
$P$ does not efficiency-dominate $Q$, and $Q$ does not 
efficiency-dominate $P$.
\end{corollary}

\proof If $P$ and $Q$ have identical eigenvalues, then
$\trace(P) = \trace(Q)$, so this follows immediately from
Theorem~\ref{tracesame}.
\qed

% The following lemma along with Propositions~\ref{eigenefprop}
% and~\ref{eigensame} will allow us to show that some reversible chains
% cannot be efficiency-dominated by \un{any} other reversible chain.

We next present a lemma about minimal values of $\trace(P)$.

\begin{lemma}\label{eigensum}
For any transition matrix $P$ on a finite state space $S$, 
for which $\pi$ is a stationary distribution,
the sum of the diagonal elements of $P$ --- that is, $\trace(P)$ ---
must be at least $\max\,(0,\ (2\pi_{\max}-1)\,/\,\pi_{\max})$, where 
$\pi_{\max}=\max_x \pi(x)$.  Furthermore, any $P$ attaining this
minimum value will have at most one non-zero value on its diagonal,
and any such non-zero diagonal value will be for a state $x^*$
for which $\pi(x^*)=\pi_{\max} > 1/2$.
\end{lemma}

\proof 
The statement is trivial when $\pi_{\max} \le 1/2$, since the lower
limit on $\trace(P)$ is then zero, and any such $P$ has all zeroes
on the diagonal.
Otherwise, if
$x^*$ is such that $\pi(x^*)=\pi_{\max} > 1/2$, then 
stationarity implies that 
\begin{eqnarray*}
  \pi_{\max} \ = \ \pi(x^*) & = & \sum_{x\in S} \pi(x) P(x,x^*)
  \ =\ \pi(x^*)P(x^*,x^*)\ +\!\!\sum_{x\in S,\,x\ne x^*}\!\!\!\pi(x)P(x,x^*)
  \\[3pt]
& \le & \pi(x^*)P(x^*,x^*)\ +\!\!\sum_{x\in S,\, x\ne x^*}\!\!\!\pi(x)
\ =\ \pi(x^*)P(x^*,x^*)\ +\ (1-\pi(x^*)) \\[3pt]
& =& \pi_{\max}\,P(x^*,x^*)\ + (1 - \pi_{\max})
\, .
\end{eqnarray*}
It follows that $P(x^*,x^*) \ge (2\pi_{\max}-1)\,/\,\pi_{\max}$,
hence $\trace(P) \ge \max\,(0,\ (2\pi_{\max}-1)\,/\,\pi_{\max})$.
Furthermore, if $\trace(P)=\max\,(0,\ (2\pi_{\max}-1)\,/\,\pi_{\max})$,
then $\trace(P)=P(x^*,x^*)$, and hence all other values on the diagonal
of $P$ must be zero.\vspace{-10pt}
\qed

\remark For any $\pi$, the minimum value of $\trace(P)$ in 
Lemma~\ref{eigensum} is attainable, and can indeed be attained
by a $P$ that is reversible.  Several methods for constructing
such a $P$ are discussed in the companion paper~\cite{radfordnew}, 
including, for example, the ``shifted tower'' method of~\cite{suwa},
which produces a reversible $P$ when the shift is by $1/2$.

\bigskip

We can now state a criterion for a reversible chain to not be
efficiency-dominated by any other reversible chain:

\begin{theorem}\label{nodom}
If $P$ is the transition matrix for an irreducible Markov chain 
on a finite state space that is reversible with respect to $\pi$, 
and the sum of the eigenvalues of $P$ (equivalently, 
the trace of $P$) equals
$\max\,(0,\ (2\pi_{\max}-1)\,/\,\pi_{\max})$, where 
$\pi_{\max}=\max_x \pi(x)$, then no other 
reversible chain can efficiency-dominate $P$.
\end{theorem}

\proof
A reversible chain, $Q$, not equal to $P$, that efficiency-dominates
$P$, must by Theorem~\ref{tracesame} have $\trace(Q) < \trace(P)$.
But by Lemma~\ref{eigensum}, $\trace(P)$ is as small as possible.
So there can be no reversible $Q$ that efficiency-dominates $P$.
\qed

\medskip

As an example of how this theorem can be applied, if the state space is
$S=\{1,2,3\}$, with $\pi(1)=\pi(2)=1/5$ and $\pi(3)=3/5$, for which
$\pi_{\max}=3/5$, then the transition matrix
$$
  P_1\ =\ \pmatrix{0 & 0 & 1 \NR 0 & 0 & 1 \NR 1/3 & 1/3 & 1/3}
$$
is reversible with respect to $\pi$, and has eigenvalues of $1, 0, -2/3$, 
which sum to 1/3 (the trace).  By Theorem~\ref{nodom}, $P_1$ cannot
be efficiency-dominated by any other reversible transition matrix, since 
its sum of eigenvalues is equal to $(2\pi_{\max}-1)\,/\,\pi_{\max}$.

On the other hand, consider the following
transition matrix, reversible with respect to the same $\pi$:
$$
  P_2\ =\ \pmatrix{0 & 1/4 & 3/4 \NR 1/4 & 0 & 3/4 \NR 1/4 & 1/4 & 1/2}
$$
$P_2$ has eigenvalues of $1, -1/4, -1/4$, which sum to $1/2$,
greater than $(2\pi_{\max}-1)\,/\,\pi_{\max}$, so 
Theorem~\ref{nodom} does not apply.  However, $P_2$ is an instance
of a transition matrix constructed according to a procedure 
of Frigessi, Hwang, and Younes \cite[Theorem~1]{frigessi}, which
they prove has the property that the transition matrix produced
has the smallest possible value for $\lambda_2$, and subject
to having that value for $\lambda_2$, the smallest possible 
value for $\lambda_3$, etc. We can therefore again conclude from
Proposition~\ref{eigenefprop} and Corollary~\ref{eigensame} that 
no other reversible chain can efficiency-dominate $P_2$.

It's easy to see that any reversible $P$ with at least two non-zero
diagonal elements, say $P(x,x)$ and $P(y,y)$, can be
efficiency-dominated by a chain, $Q$, that is the same as $P$ except
that these diagonal elements are reduced, allowing $Q(x,y)$ and $Q(y,x)$
to be greater than $P(x,y)$ and $P(y,x)$, so that $Q$ Peskun-dominates
$P$.  Theorem~\ref{nodom} shows that \un{some} reversible $P$ in which
only a single diagonal element is non-zero cannot be
efficiency-dominated by any other reversible chain.
We know of no examples of a reversible $P$ with only one non-zero diagonal
element that is dominated by another reversible chain, but we do not have
a proof that this is impossible.
This leads to:

\newpartitle{Open Problem}
Does there exists a reversible $P$ with only one non-zero
diagonal element, which is efficiency-dominated by some other
reversible chain?

\section{Re-deriving Peskun's Theorem}\label{sec-peskun}

Recall that $P$ \un{Peskun-dominates} $Q$ if
$P(x,y) \ge Q(x,y)$ for all $x\not=y$ ---
i.e., that $Q\!-\!P$ has all non-positive entries off the diagonal
(and hence also that $Q-P$ has all non-negative entries on the diagonal).
It is known through several complicated 
proofs~\cite{peskun, tierney2, radfordnonrev}
that if $P$ Peskun-dominates $Q$, then
$P$ efficiency-dominates $Q$.
We will see here that once Theorem~\ref{eigenthm} has been established, this
fact can be shown easily.

\begin{proposition}\label{peskunprop}
If $P$ and $Q$ are irreducible, and both are reversible with respect to 
some $\pi$, and $P$ Peskun-dominates $Q$, then $P$ efficiency-dominates $Q$.
\end{proposition}

To prove Proposition~\ref{peskunprop}, we begin with
a simple eigenvalue lemma.

\begin{lemma}\label{rowsumlemma}
If $Z$ is an $n \times n$ matrix with $z_{ii} \ge 0$
and $z_{ij} \le 0$ for all $i\not=j$, and
row-sums $\sum_j z_{ij} = 0$ for all $i$,
then all eigenvalues of $Z$ must be non-negative.
\end{lemma}

\proof
Suppose $Zv=\lambda v$.  Find the index $j$ which maximizes $|v_j|$,
i.e.\ such that $|v_j| \ge |v_k|$ for all $k$.
We can assume $v_j > 0$ (if not, replace $v$ by $-v$),
so $v_j \ge |v_k|$ for all $k$.
Then
$$
    \lambda \, v_j \ = \ (Zv)_j \ = \ \sum_i z_{ji} v_i
		\ = \ z_{jj} v_j + \sum_{i \not= j} z_{ji} v_i
                \ \ge \ z_{jj} v_j - \sum_{i \not= j} |z_{ji}| \, |v_i|
$$
$$
                \ \ge \ z_{jj} v_j - \sum_{i \not= j} |z_{ji}| \, v_j
                % \ = \ v_j \Big( z_{jj} - \sum_{i \not= j} |z_{ji}| \Big)
                \ = \ v_j \Big( z_{jj} + \sum_{i \not= j} z_{ji} \Big)
                \ = \ v_j (0)
                \ = \ 0
\, .
$$
So, $\lambda \, v_j \ge 0$.
Hence, since $v_j>0$, we must have $\lambda \ge 0$.
\qed\vspace{-5pt}

\newpartitle{Proof of Proposition~\ref{peskunprop}}
Let $Z=Q\!-\!P$.  Since $P$ Peskun-dominates $Q$,
$z_{ii} = Q(i,i) - P(i,i) \ge 0$
and $z_{ij} = Q(i,j) - P(i,j) \le 0$ for all $i\not=j$.
Also $\sum_j z_{ij} = \sum_j P(i,j) - \sum_j Q(i,j) = 1 - 1 = 0$.
Hence, by
Lemma~\ref{rowsumlemma}, $Z=Q\!-\!P$ has all eigenvalues non-negative.
Hence, by Theorem~\ref{eigenthm}, $P$ efficiency-dominates $Q$.
\qed

\remark Proposition~\ref{peskunprop} can also be proven by transforming
$Q$ into $P$ one step at a time, in the sequence $Q,Q',Q'',\ldots,P$,
with each matrix in the sequence efficiency-dominating the previous matrix.
At each step, say from $Q'$ to $Q''$, two of the off-diagonal transition 
probabilities that differ between $Q$ and $P$, say
those involving states $x$ and $y$, will be increased from $Q(x,y)$ to
$P(x,y)$ and from $Q(y,x)$ to $P(y,x)$, while $Q''(x,x)$ and $Q''(y,y)$ will 
decrease compared to $Q'(x,x)$ and $Q'(y,y)$.  The difference $Q'-Q''$
will be zero except for a $2 \times 2$ submatrix involving states $x$
and $y$, which will have the form 
$\Big(\begin{array}{rr} a & \!-a \\[-2pt] \!\!-b & b \end{array}\Big)$ 
for some $a,b>0$, which has non-negative eigenvalues of 0 and $a+b$.
Hence, by Theorem~\ref{eigenthm}, $Q''$ efficiency-dominates $Q'$.  
Since this will be true for all the steps from $Q$ to $P$, 
transitivity (see Theorem~\ref{partialord}(c)) implies that
$P$ efficiency-dominates $Q$. \qed

Note that the converse to Proposition~\ref{peskunprop}
is \un{false}.  For example, let $S=\{1,2,3\}$, and
$$
P = \pmatrix{ 0 & 1/2 & 1/2 \NR 1 & 0 & 0 \NR 1 & 0 & 0 }
, \quad
Q = \pmatrix{ 1/2 & 1/4 & 1/4 \NR 1/2 & 1/4 & 1/4 \NR 1/2 & 1/4 & 1/4 }
, \quad
Q\!-\!P =
\pmatrix{ +1/2 & -1/4 & -1/4 \NR -1/2 & +1/4 & +1/4 \NR -1/2 & +1/4 & +1/4 }
\, .
$$
Here, $P$ does \un{not} Peskun-dominate $Q$, since, for example, 
$Q(2,3)=1/4 > 0 = P(2,3)$. 
However, the eigenvalues of $Q\!-\!P$ are $1,0,0$, 
all of which are non-negative, so $P$ \un{does} efficiency-dominate $Q$.
Furthermore, Theorem~\ref{partialord}(b) implies
that $P$ is strictly better than $Q$ --- there is some $f$ for which 
$v(f,P)<v(f,Q)$. (For example, the indicator function for the first state,
which has asymptotic variance zero using $P$, and asymptotic variance $1/4$
using $Q$.)
Peskun dominance therefore does not capture all instances of efficiency 
dominance that we are interested in, which motivates our investigation here.

One should note, however, that all our results concern only reversible
chains.  Non-reversible chains are often used, either in a deliberate
attempt to improve performance (see \cite{radfordnonrev}), or somewhat
accidentally, as a result of combining methods sequentially rather
than by random selection.
Extensions of Peskun ordering to non-reversible chains are
considered in \cite{AndrieuLivingstone}.
% Investigation of non-reversible chains is another area for future work.

\section{Elementary Proof of Lemma~\ref{effequivlemma}}
\label{sec-effequivproof}

We conclude by presenting the promised elementary proof of 
Lemma~\ref{effequivlemma} above, which states the surprising fact that,
for any reversible irreducible transition matrices $P$ and $Q$,
$\inn{f}{Pf} \le \inn{f}{Qf}$ for all $f$
\un{if and only if}
$\inn{f}{P(I\!-\!P)^{-1}f} \le \inn{f}{Q(I\!-\!Q)^{-1}f}$ for all $f$.

As observed in \cite[pp.~16--17]{mirageyer}, this proposition follows
from the more general result
of Bendat and Sherman \cite[p.~60]{bendat},
using results of L\"owner \cite{lowner},
which states that that if
$h(x) = {ax+b \over cx+d}$ where $ad-bc>0$,
and
$J$ and $K$ are any two self-adjoint operators
with spectrum contained in
$(-\infty,-d/c)$
or in $(-d/c,\infty)$,
then if
$\inn{f}{Jf} \le \inn{f}{Kf}$ for all $f$, then also
$\inn{f}{h(J)f} \le \inn{f}{h(K)f}$ for all $f$.
In particular, choosing $a=d=1$, $b=0$, and $c=-1$ gives that
$h(J) = {J \over I-J}$, so
if $\inn{f}{Jf} \le \inn{f}{Kf}$ for all $f$ then
$\inn{f}{{J \over I-J} \, f} \le \inn{f}{{K \over I-K} \, f}$ for all $f$.
Conversely, choosing $a=c=d=1$ and $b=0$ gives
that $h({J \over I-J}) = J$, so if
$\inn{f}{{J \over I-J} \, f} \le \inn{f}{{K \over I-K} \, f}$
for all $f$
then $\inn{f}{Jf} \le \inn{f}{Kf}$ for all $f$,
finishing the proof.

However, the proof in \cite{bendat} is very technical,
requiring analytic continuations of transition functions into
the complex plane.  Instead, we now present an elementary proof
of Lemma~\ref{effequivlemma}.
(See also \cite[Chapter~V]{Bhatia}.)
We begin with some lemmas
about operators on a finite vector space $\V$, e.g.\ $\V=L^2_0(\pi)$.

\begin{lemma}\label{Zlemma}
If $X,Y,Z$ are
% $n \times n$ matrices,
operators on a finite vector space $\V$,
with $Z$ self-adjoint,
and $\inn{f}{Xf} \le \inn{f}{Yf}$ for all $f\in\V$,
then $\inn{f}{ZXZf} \le \inn{f}{ZYZf}$ for all $f\in\V$.
\end{lemma}

\proof
Since $Z$ is self-adjoint, making the substitution $g=Zf$ gives
$$
\inn{f}{ZXZf}
\ = \ \inn{Zf}{XZf}
\ = \ \inn{g}{Xg}
\ \le \ \inn{g}{Yg}
\ = \ \inn{Zf}{YZf}
\ = \ \inn{f}{ZYZf}
\, .
\eqqed
$$

\medskip

Next, say a self-adjoint matrix
$J$ is \un{strictly positive} if
$\inn{f}{Jf} > 0$ for all non-zero $f\in\V$.
Since $\inn{f}{Jf} = \sum_i (a_i)^2 \lambda_i$
(see~Section~\ref{sec-prelim}),
this is equivalent
to $J$ having all eigenvalues positive.

\begin{lemma}\label{inverselemmaoneway}
If $J$ and $K$ are strictly positive self-adjoint
% $n\times n$ matrices,
operators on a finite vector space $\V$,
% then
and
$\inn{f}{Jf} \le \inn{f}{Kf}$ for all $f\in\V$,
% \un{if and only if}
then $\inn{f}{J^{-1}f} \ge \inn{f}{K^{-1}f}$ for all $f\in\V$.
\end{lemma}

\proof
% We first show the forward implication, that 
% $\inn{f}{Jf} \le \inn{f}{Kf}$ for all $f\in\V$ implies that
% $\inn{f}{J^{-1}f} \ge \inn{f}{K^{-1}f}$ for all $f\in\V$.
Note (see Section~\ref{sec-prelim}) that $K^{-1/2}$ and 
$K^{-1/2}JK^{-1/2}$ are self-adjoint, and have
eigenvalues that are positive, since $J$ and $K$ are 
strictly positive.
Applying Lemma~\ref{Zlemma} with $X=J$, $Y=K$, and $Z=K^{-1/2}$
then gives that for all $f\in\V$,
$$
\inn{f}{K^{-1/2}JK^{-1/2}f}
\ \le \ \inn{f}{K^{-1/2}KK^{-1/2}f}
\ = \ \inn{f}{If}
\ = \ \inn{f}{f}
\, .
$$
It follows that all the eigenvalues of $K^{-1/2}JK^{-1/2}$ are 
in $(0,1]$ (since $\inn{v}{Av}\le \inn{v}{v}$ and $Av=\lambda v$ 
with $v\ne0$ imply $\lambda\le 1$).
Hence, its inverse
$\big( K^{-1/2}JK^{-1/2} \big)^{-1}$ has eigenvalues all $\ge 1$, 
so
$
\inn{f}{\big( K^{-1/2}JK^{-1/2} \big)^{-1}f}
\ \ge \
\inn{f}{If}
\, ,
$
and therefore
$$
\inn{f}{If}
\ \le \ \inn{f}{\big( K^{-1/2}JK^{-1/2} \big)^{-1}f}
\ = \ \inn{f}{K^{1/2}J^{-1}K^{1/2}f}
\, .
$$
Then, applying Lemma~\ref{Zlemma} again
with $X=I$, $Y=K^{1/2}J^{-1}K^{1/2}$, and $Z=K^{-1/2}$ gives
$$
\inn{f}{K^{-1/2}IK^{-1/2}f}
\ \le \ \inn{f}{K^{-1/2}\, \big(K^{1/2}J^{-1}K^{1/2}\big)\, K^{-1/2}f}
\, .
$$
That is,
$
\inn{f}{K^{-1}f} \ \le \ \inn{f}{J^{-1}f}\ \ \mbox{for all $f \in \V$,}
$
giving the result.
\qed

\def\fifth{{1 \over 5}}

We emphasise that the equivalence in
Lemma~\ref{inverselemmaoneway}
is only
for \un{all} $f$ at once, not for individual $f$.
For example, if
$J=I$ and $K=\diag(5,\fifth)$ and $f=(1,1)$, then
$\inn{f}{Jf} = 2 \le 5 + \fifth = \inn{f}{Kf}$,
but $\inn{f}{J^{-1}f} = 2 \not\ge \fifth + 5 = \inn{f}{K^{-1}f}$.
This illustrates why the proofs of
Lemma~\ref{inverselemmaoneway}
and Lemma~\ref{effequivlemma} are not as straightforward as one might think.

\remark
Lemma~\ref{inverselemmaoneway} can be partially proven more directly.
If $f=\sum_i a_i v_i$,
% then since $\phi(x)=1/x$ is \un{convex} on $\{x>0\}$
% (because $\phi''(x) = 2x^{-3} >0$ there),
\un{Jensen's Inequality}
% with $\phi(x)=1/x$
gives
$
% \phi\Big( \sum_i (a_i)^2 \lambda_i \Big)
\Big( \sum_i (a_i)^2 \lambda_i \Big)^{-1}
\le
% \sum_i (a_i)^2 \phi(\lambda_i)
\sum_i (a_i)^2 \, (\lambda_i)^{-1}
$,
so we always have
$
% {1 \over \inn{f}{Kf}}
% \ \le \
{1 \over \inn{f}{Jf}}
\le \inn{f}{J^{-1}f}
\, .
$
Hence, if $\inn{f}{Jf} \le \inn{f}{Kf}$, then
% ${1 \over \inn{f}{Kf}} \ \le \ \inn{f}{J^{-1}f}$.
$\inn{f}{J^{-1}f} \ge {1 \over \inn{f}{Kf}}$.
If $f$ is an \un{eigenvector} of $K$, then
$
\inn{f}{K^{-1}f}
\ = \ {1 \over \inn{f}{Kf}}
$,
% This provides a more direct proof of
so this shows directly that
% $\inn{f}{K^{-1}f} \le \inn{f}{J^{-1}f}$.
$\inn{f}{J^{-1}f} \ge \inn{f}{K^{-1}f}$.
% Lemma~\ref{inverselemmaoneway}.
% in the case where $f$ is an eigenvector of $K$.
However, it is unclear how to extend this argument to other $f$.
\qed

Applying Lemma~\ref{inverselemmaoneway} twice gives a
(stronger) two-way equivalence:

\begin{lemma}\label{inverselemma}
If $J$ and $K$ are strictly positive self-adjoint
% $n\times n$ matrices,
operators on a finite vector space $\V$, then
$\inn{f}{Jf} \le \inn{f}{Kf}$ for all $f\in\V$
\un{if and only if} $\inn{f}{J^{-1}f} \ge \inn{f}{K^{-1}f}$ for all $f\in\V$.
\end{lemma}

% that $\inn{f}{J^{-1}f} \ge \inn{f}{K^{-1}f}$ for all $f\in\V$ implies
% $\inn{f}{Jf} \le \inn{f}{Kf}$ for all $f\in\V$,

\proof
The forward implication is
Lemma~\ref{inverselemmaoneway}.
And, the reverse implication follows from
Lemma~\ref{inverselemmaoneway}
by replacing $J$
with $K^{-1}$ and replacing $K$ with $J^{-1}$.
\qed

%\smallskip

\iffalse

\medskip Applying Lemma~\ref{inverselemma}
again with $J$ replaced by $K^{-1}$
and $K$ replaced by $J^{-1}$
gives:

\begin{lemma}\label{equivlemma}
If $J$ and $K$ are strictly positive self-adjoint
% $n\times n$ matrices,
operators on a finite vector space $\V$,
then $\inn{f}{Jf} \le \inn{f}{Kf}$ for all $f$
\un{if and only if}
$\inn{f}{J^{-1}f} \ge \inn{f}{K^{-1}f}$ for all $f$.
\end{lemma}

\fi

\smallskip Using Lemma~\ref{inverselemma}, we easily obtain:

\newpartitle{Proof of Lemma~\ref{effequivlemma}}
Recall that we can restrict to $f\in L^2_0(\pi)$, so $\pi(f)=0$,
and $f$ is orthogonal to the eigenvector corresponding to
eigenvalue~1.  On that restricted subspace,
the eigenvalues of $P$ and $Q$ are contained in $[-1,1)$.
Hence, the eigenvalues of $I\!-\!P$ and $I\!-\!Q$ are contained in
$(0,2]$, and in particular are all strictly positive.
So, $I\!-\!P$ and $I\!-\!Q$ are strictly positive self-adjoint operators.

Now,
$\inn{f}{Pf} \le \inn{f}{Qf}$ for all $f \in L^2_0(\pi)$ is 
equivalent to 
$$\inn{f}{(I\!-\!P)f} = \inn{f}{f}-\inn{f}{Pf}
\ \ge\ \inn{f}{f}-\inn{f}{Qf} = \inn{f}{(I\!-\!Q)f}\ \ 
\mbox{for all $f \in L^2_0(\pi)$}.$$
Then, by
Lemma~\ref{inverselemma} with $J = I\!-\!Q$ and $K = I\!-\!P$, this
is equivalent to
$$\inn{f}{(I\!-\!P)^{-1}f}\ \le\ \inn{f}{(I\!-\!Q)^{-1}f}\ \ 
\mbox{for all $f \in L^2_0(\pi)$}.$$
Since $(I\!-\!P)^{-1} = P(I\!-\!P)^{-1} + (I-P)(I-P)^{-1}
= P(I\!-\!P)^{-1} + I$ and similarly $(I\!-\!Q)^{-1} = Q(I\!-\!Q)^{-1} + I$,
this latter is equivalent to
$$\inn{f}{P(I\!-\!P)^{-1}f}\ \le\ \inn{f}{Q(I\!-\!Q)^{-1}f}\ \ 
\mbox{for all $f \in L^2_0(\pi)$},$$
which completes the proof.
\qed

\medskip
\ackn
We are very grateful to
Heydar Radjavi for the proof of Lemma~\ref{inverselemmaoneway} herein,
and thank Rajendra Bhatia, Gareth Roberts, Daniel Rosenthal,
and Peter Rosenthal for several very
helpful discussions about eigenvalues.
We thank the anonymous referee for a number of useful suggestions
which helped to improve the paper.

\end{document}